\documentclass[a4paper,10pt]{article}
\usepackage{mathtext}
\usepackage[latin2]{inputenc}
\usepackage[english]{babel}
\usepackage{amsmath}
\usepackage{amsfonts}
\usepackage{amssymb}
\usepackage{mathrsfs}
\usepackage{amsthm}
\usepackage{euscript}
\usepackage{graphicx}
\usepackage{float}
\DeclareMathOperator{\RR}{\mathbb{R}}
\DeclareMathOperator{\TT}{\mathbb{T}}

\newcommand{\eLd}[1]{\|{#1}\|_{L_1(\mathbb{T}^d)}}
\newcommand{\eL}[1]{\left\|{#1}\right\|_{L_1(\mathbb{T}^2)}}

\renewcommand{\leq}{\leqslant}
\renewcommand{\geq}{\geqslant}

\theoremstyle{plain}\newtheorem{T1}{Theorem}
\theoremstyle{plain}
\theoremstyle{plain}\newtheorem{L1}{Lemma}
\theoremstyle{plain}
\theoremstyle{plain}
\theoremstyle{plain}
\theoremstyle{plain}

\begin{document}
\author{Krystian Kazaniecki \and Michal Wojciechowski}
\title{Ornstein's non-inequalities Riesz product approach}
\maketitle
\begin{abstract}We provide a new technique to prove Ornstein's non-inequalities for derivatives with some geometrical dependence on their indexes.

\end{abstract}

\vskip 3mm
In \cite{MR0149331} D. Ornstein studied behavior of norms of the homogeneous differential operators with the same degree of homogeneity. He proved lack of a priori estimates for linearly independent differential operators.
We show alternative way of proving Orstein's non-inequalities in some special cases. 
\begin{T1}
 Assume $\alpha_0, \ldots, \alpha_n$ are multindexes in $\RR^d$. If there exists a pair of vectors $\Gamma,\; \Lambda\in(\mathbb{N}\cup\{0\})^d$ for which following occurs
 \[
 \langle \alpha_0,\, \Lambda \rangle= \langle \alpha_1,\, \Lambda \rangle=\ldots= \langle \alpha_m,\, \Lambda \rangle,
 \]
and
\[
\langle \alpha_0,\, \Gamma \rangle< \langle \alpha_1,\, \Gamma \rangle<\langle \alpha_2,\, \Gamma \rangle\leq\ldots\leq\langle \alpha_m,\, \Gamma \rangle,
\]
Then for every $K>0$ there exists $f\in C^\infty(\TT^d)$ such that 
\begin{equation*}
\eLd{D^{\alpha_1} f}\geq K \sum_{j\neq 1} \eLd{D^{\alpha_j}f}.
\end{equation*}
\end{T1}
Instead of giving full proof, we limit ourself to the special, yet representative, case. We prove that for every $K>0$ there exists $f\in C^\infty(\TT^d)$ such that 
\begin{equation}\label{nierownosc}
\eL{\frac{\partial^5}{\partial x_1^3 \partial x_2^2} f}\!\!\geq K \left(\eL{\frac{\partial^4}{\partial x_1^4 }f}\!+\! \eL{\frac{\partial^6}{\partial x_1^2\partial x_2^4 }f}\!+\!\eL{\frac{\partial^7}{\partial x_1 \partial x_2^6} f}\!+\! \eL{\frac{\partial^8}{\partial x_2^8} f}\right).
\end{equation}
\begin{proof}
 We fix $K>0$ and $n> 64 K^2 C^{-2}$. We will construct trigonometric polynomial, whose one of the derivatives behaves like the modified Riesz product
 \begin{equation*}
R_n(x)=-1 + \Pi_{k=1}^{n}\left(1+\cos(2\pi \langle x \, , \, a_k\rangle \right) ,
 \end{equation*}
where $a_k\in\mathbb{Z}^2$. By induction we can choose $a_k\in\mathbb{Z}^2$ for $k=1,\ldots, n$ such that  
\begin{equation}\label{wlasnosci}
\begin{aligned}
 \left|\left(\frac{\left( a_k(2)+ \sum_{j=1} \epsilon_j a_{j}(2)\right)^2}{a_k(1)+ \sum_{j=1}\epsilon_j a_{j}(1)}\right)^l - \left(\frac{\sigma_k}{\sqrt{n}}\right)^l\right|&\leq \frac{1}{3^n},\\
 \|a_{k+1}|&> M_n |a_{k}| ,\\
 a_k(1)+ \sum_{j=1}\epsilon_j a_{j}(1)&\neq 0,
 \end{aligned}
\end{equation}
for every $k\in\{1,2,\ldots,n\}$, $l\in\{1,2,\ldots,m\}$ and $\epsilon_j\in \{-1,0,1\}$, where $n> (8K+1)^2 C^{-2}$ and the values of parameters $M_n$, $C$ and $\sigma_j$ are determined by the following lemma. 
\begin{L1}\label{lematpolostateczny}(cf. \cite{MR1649869})
 There is $C>0$ such that for every $m\in\mathbb{N}^{+}$ there exists $M=M(m)$ and a sequence $\{\sigma_j\}_{j=1}^{m}\in \{0,1\}^{m}$ such that
\begin{equation*}
\left\|\sum_{j=1}^{m} \sigma_j\cos \left(2\pi \langle d_j,\,\xi \rangle\right)\prod_{1\leq k< j}\left(1+\cos\left(2\pi\langle d_{k},\, \xi\rangle\right)\right)\right\|_{L_1(\mathbb{T}^d)}\geq C n
\end{equation*}
 whenever $\{d_k\}_{k=1}^{m}\subset\mathbb{Z}^d$ satisfies $|d_{k+1}|> M_m |d_{k}|$ for $k=1,\ldots,m-1$.
\end{L1} 
This inequality was generalized by R. Lata\l a in \cite{info}. We define family of sets 
\begin{equation*}
A_k=\{q : q=a_k + \sum_{j=1}^{k-1} \epsilon_j a_j,\; \epsilon_j=-1,0,1\;\}. 
\end{equation*}
For $q=\sum_{j=1}^{n} \epsilon_j(q) a_j$ we put $r(q)=\# \{j:  \epsilon_j\neq 0\}$.
Let $Z$ be the polynomial given by the formula
\begin{equation*}
Z(x)=\sum_{k=1}^{n}\sum_{q\in A_k\cup-A_k} \frac{1}{q(1)^4}\frac{1}{2^{r(q)}} e^{2i\pi\langle q,\, x\rangle}.
\end{equation*}
Simple calculation gives 
\begin{equation*}
D^{\alpha_0} Z(x)=\sum_{k=1}^{n}\sum_{q\in A_k\cup-A_k} \frac{1}{2^{r(q)}} e^{2i\pi\langle q,\, x\rangle}= R_n(x).
\end{equation*}
Hence 
\begin{equation*}
\eL{D^{\alpha_0} Z_n (x)}\leq 2
\end{equation*}
Since $D^{\alpha_l} Z(x)=D^{\alpha_l - \alpha_0} D^{\alpha_0}Z(x)=D^{l(\alpha_1 - \alpha_0)}R_n$, for $m\geq 1$ we have
\begin{equation*}
\begin{split}
i^{-l}D^{\alpha_l} Z(x)=\sum_{k=1}^{n}\sum_{q\in A_k\cup-A_k} \left(\frac{q(2)^{2}}{q(1)}\right)^l\frac{1}{2^{r(q)}} e^{2i\pi\langle q,\, x\rangle}= I_l + I\!I_l,
\end{split}
\end{equation*}
where 
\begin{equation*}
\begin{split}
I_l &= \sum_{k=1}^{n}\sum_{q\in A_k\cup-A_k} \left(\frac{q(2)^{2l}}{q(1)^l}- \left(\frac{\sigma_k}{\sqrt{n}}\right)^l\right)\frac{1}{2^{r(q)}} e^{2i\pi\langle q,\, x\rangle}, \\
I\!I_l&= \sum_{k=1}^{n}\sum_{q\in A_k\cup-A_k} \left( \frac{\sigma_k}{\sqrt{n}}\right)^l\frac{1}{2^{r(q)}} e^{2i\pi\langle q,\, x\rangle}.
\end{split}
\end{equation*}
Since \eqref{wlasnosci} is satisfied and $\# A_k = 3^{k-1}$,
\begin{equation*}
 \eL{I_l}\leq 1
\end{equation*}
It is easy to check that 
\begin{equation*}
I\!I_l= n^{-\frac{l}{2}} \sum_{k=1}^{n} \sigma_k^l\cos \left(2\pi \langle a_k,\,\xi \rangle\right)\prod_{1\leq j< k}\left(1+\cos\left(2\pi\langle a_{j},\, \xi\rangle\right)\right).
\end{equation*}
By the triangle inequality for $m\geq 2$ we have
\begin{equation*}
\eL{I\!I_m}\leq n^{-\frac{m}{2}} n \leq 1
\end{equation*}
For $m=1$ by Lemma \ref{lematpolostateczny}, we get
\begin{equation*}
 \eL{ I\!I_1 }\geq Cn n^{-\frac{1}{2}} = C\sqrt{n}> 8K+1,
\end{equation*}
which proves \eqref{nierownosc}.

\end{proof}
\bibliographystyle{plain}
\bibliography{cont}

\begin{thebibliography}{1}

\bibitem{info}
Rafa{\l} Lata{\l}a.
\newblock {$L_1$-norm of combinations of products of independent random
  variables}.
\newblock {\em Israel J. Math.}, to appear.

\bibitem{MR0149331}
Donald Ornstein.
\newblock {A non-equality for differential operators in the {$L_{1}$} norm.}
\newblock {\em Arch. Rational Mech. Anal.}, 11:40--49, 1962.

\bibitem{MR1649869}
Micha{\l} Wojciechowski.
\newblock {On the strong type multiplier norms of rational functions in several
  variables}.
\newblock {\em Illinois J. Math.}, 42(4):582--600, 1998.

\end{thebibliography}

\end{document}